\documentstyle[12pt]{article}
\begin{document}
\newpage
\pagestyle{empty}
\setcounter{page}{0}
\renewcommand{\theequation}{\thesection.\arabic{equation}}
\newcommand{\sect}[1]{\setcounter{equation}{0}\section{#1}}
\newfont{\twelvemsb}{msbm10 scaled\magstep1}
\newfont{\eightmsb}{msbm8}
\newfont{\sixmsb}{msbm6}
\newfam\msbfam
\textfont\msbfam=\twelvemsb
\scriptfont\msbfam=\eightmsb
\scriptscriptfont\msbfam=\sixmsb
\catcode`\@=11
\def\Bbb{\ifmmode\let\next\Bbb@\else
  \def\next{\errmessage{Use \string\Bbb\space only in math mode}}\fi\next}
\def\Bbb@#1{{\Bbb@@{#1}}}
\def\Bbb@@#1{\fam\msbfam#1}
\newfont{\twelvegoth}{eufm10 scaled\magstep1}
\newfont{\tengoth}{eufm10}
\newfont{\eightgoth}{eufm8}
\newfont{\sixgoth}{eufm6}
\newfam\gothfam
\textfont\gothfam=\twelvegoth
\scriptfont\gothfam=\eightgoth
\scriptscriptfont\gothfam=\sixgoth
\def\frak{\frak@}
\def\frak@#1{{\fam\gothfam{{#1}}}}
\def\frak@@#1{\fam\gothfam#1}
\catcode`@=12
\def\arcsinh{\mathop{\rm arcsinh}\nolimits}
\def\CC{{\Bbb C}}
\def\NN{{\Bbb N}}
\def\QQ{{\Bbb Q}}
\def\RR{{\Bbb R}}
\def\ZZ{{\Bbb Z}}
\def\cA{{\cal A}}          \def\cB{{\cal B}}          \def\cC{{\cal C}}
\def\cD{{\cal D}}          \def\cE{{\cal E}}          \def\cF{{\cal F}}
\def\cG{{\cal G}}          \def\cH{{\cal H}}          \def\cI{{\cal I}}
\def\cJ{{\cal J}}          \def\cK{{\cal K}}          \def\cL{{\cal L}} 
\def\cM{{\cal M}}          \def\cN{{\cal N}}          \def\cO{{\cal O}}
\def\cP{{\cal P}}          \def\cQ{{\cal Q}}          \def\cR{{\cal R}} 
\def\cS{{\cal S}}          \def\cT{{\cal T}}          \def\cU{{\cal U}}
\def\cV{{\cal V}}          \def\cW{{\cal W}}          \def\cX{{\cal X}}
\def\cY{{\cal Y}}          \def\cZ{{\cal Z}}
\def\id{\mbox{id}}
\def\ggo{{\frak g}_{\bar 0}}
\def\uqggo{\cU_q({\frak g}_{\bar 0})}
\def\uqggp{\cU_q({\frak g}_+)}
\def\typeA{{\em type $\cA$}}
\def\typeB{{\em type $\cB$}}

$$
\;
$$

\vskip 2cm

\begin{center}

  {\LARGE {\bf {\sf On Nonlinear Angular Momentum Theories, Their 
		Representations and Associated Hopf Structures }}} \\[1cm]

\smallskip

{\large B. Abdesselam$^{\dagger,2}$,
  J. Beckers$^{\ddagger,3}$, A. Chakrabarti$^{\dagger,4}$
  and N. Debergh$^{\ddagger,3,5}$}

\vfill

{\em $^{\dagger}$
  Centre de Physique Th{\'e}orique
  \footnote{Laboratoire Propre du CNRS UPR A.0014  

  \indent
  $\;$$^2$abdess@orphee.polytechnique.fr,

  \indent
  $\;$$^3$beckers@vm1.ulg.ac.be,

  \indent
  $\;$$^4$chakra@orphee.polytechnique.fr,

  \indent
  $\;$$^5$Chercheur Institut Interuniversitaire des Sciences 
  Nucl{\'e}aires (Brussels).}, Ecole Polytechnique, 
  91128 Palaiseau Cedex,
  France.
  \indent
  \vskip 0.2cm
  $^{\ddagger}$ Physique Th{\'e}orique et Math{\'e}matique, Institut de 
  Physique, B5, Universit{\'e} de Li{\`e}ge, B-4000-Li{\` e}ge 1, Belgium.}

\end{center}

\vfill

\begin{abstract}

  Nonlinear $sl(2)$ algebras subtending generalized angular momentum 
  theories are studied in terms of undeformed generators and bases. We 
  construct their unitary irreducible representations in such a general 
  context. The linear $sl(2)$-case as well as its $q$-deformation are 
  easily recovered as specific examples. Two other physically interesting 
  applications corresponding to the so-called Higgs and quadratic algebras 
  are also considered. We show that these two nonlinear algebras can be 
  equipped with a Hopf structure.      
\end{abstract}

\vfill 

\rightline{CPTH-S 424-1295}
\rightline{q-alg/9601001}
\rightline{PTM-95/14 U.LIEGE}
\rightline{\sl January 96} 

\vfill 
\vfill

\newpage
\pagestyle{plain}


\sect{Introduction\label{sect:introduction}}

Quantum groups \cite{DJ} evidently appear as algebras with 
an infinite set of products of generators in the right-hand side of their 
commutation relations. If we limit the order of such products, we also get 
particular generalizations of ordinary Lie algebras that here we simply refer 
to as nonlinear algebras defined in following section: let us mention 
in particular that 
${\cal W}$-finite algebras \cite{Barbarin} belong to that category but
also that there are known examples like the Higgs algebra 
\cite{Higgs} (containing {\sl cubic} terms) and like the so-called 
{\sl quadratically} nonlinear algebras \cite{Schoutens}. Such specific
nonlinear algebras have recently been investigated by Ro${\tilde c}$ek \cite{Rocek} 
and related by Quesne \cite{Quesne} to generalized deformed 
parafermions \cite{Ohnuki} which can be exploited in the study of the 
spectra of Morse and modified P$\ddot o$schl-Teller Hamiltonians 
\cite{Daskaloyannis} as well as of parasupersymmetric Hamiltonians 
\cite{Beckers}.

The (so important) angular momentum theory being subtended by the real forms 
of the complex Lie (Cartan) algebra $A_{1}$ \cite{Cornwell}, we are here 
interested in some generalizations of this angular momentum theory to the 
nonlinear extensions of $ A_{1}$. In particular, we plan to study the 
representations associated with such nonlinear algebras. This is the 
first purpose of our study. The second one is connected with the possibility 
of endowing where it is possible these nonlinear algebras with a 
Hopf structure 
\cite{DJ}. Consequently, the contents are distributed as follows.   

In {\sl section} \ref{sect:representation}, we study a specific series 
(admitting only odd powers) of {\sl nonlinear} $sl(2)$ algebras 
subtending generalized angular momentum theories and construct their 
unitary irreducible representations. In {\sl section} 
\ref{sect:generalisation}, 
we give a generalization of nonlinear algebras when the starting point was
${\cal U}_{q}(sl(2))$. In {\sl section} \ref{sect:sl2q-case}, we show that 
the linear $sl(2)$-case as well as its $q$-deformation 
\cite{DJ} are particular examples of our developments. Some comments about
the Hopf structure of thess nonlinear algebras is given in section 
\ref{sect:Hopf-structure}. The specific {\sl cubic} context and some comments 
about the Hopf structure are then considered in {\sl section} 
\ref{sect:cubic-case}. We also show that there exist other new families 
of representations when a specific choice of the diagonal generator is 
considered. Then, we study the {\sl quadratic} context in {\sl section} 
\ref{sect:quadratic-case} by exploiting the above choice although this 
nonlinear algebra does not belong to the specific series. By the way, some 
comments about this 
quadratic $sl(2)$-algebra are also given. Finally, {\sl section} 
\ref{sect:conclusions} is devoted to general comments and conclusions in 
connection with other recent proposals.

\sect{Representation Theory of Nonlinear $sl(2)$ Algebras 
\label{sect:representation}}

In terms of the ladder generators $J_{\pm}$ and the diagonal 
one $J_{3}$, the very well-known linear $sl(2)$ algebra is characterized 
by the commutation relations \cite{Edmonds}
\begin{eqnarray}
  && [J_{+} , J_{-} ] = 2 J_{3}, \label{eq:def-sl2-1}\\
  &&  \nonumber \\  
  && [J_{3} , J_{\pm} ] = \pm J_{\pm}, \label{eq:def-sl2-2}
\end{eqnarray}
and by the Casimir operator 
\begin{equation}
{\cal C}={1\over 2}(J_{+} J_{-} + J_{-} J_{+} ) + J_{3}^{2},
\label{eq:undeformed-Casimir}
\end{equation}
acting on an orthogonal basis denoted as usual by $\lbrace |j, m\rangle 
\rbrace$. In fact, we have the well known results 
\begin{eqnarray} 
&& {\cal C} |j, m\rangle= j(j+1)|j, m\rangle,\qquad\qquad 
   j=0,\;{1 \over 2},\;1,\;{3\over 2}, \cdots \\
&& \nonumber \\  
&& J_{3}|j, m\rangle= m |j, m\rangle,\qquad\qquad 
   m=-j,\;-j+1,\cdots, j-1,\;j \\
&& \nonumber \\  
&& J_{\pm}|j, m\rangle= \sqrt{(j\mp m)(j\pm m+1)}|j, m\pm 1\rangle 
\label{eq:representation} 
\end{eqnarray}
which characterize all the unitary irreducible representations of this 
simple Lie algebra. 

Let us consider the algebras that we decide to call nonlinear $sl(2)$ 
algebras due to the nonlinear terms appearing in the right-hand side of the 
following commutation relations (in correspondence with the ones given 
by eqs. (\ref{eq:def-sl2-1}) and (\ref{eq:def-sl2-2})), i.e.
\begin{eqnarray}
&& [{\hat J}_{+}, {\hat J}_{-}]=
   \sum_{p=0}^{N}\beta_{p}(2\;{\hat J}_{3})^{2p+1}, 
\label{eq:deformed-sl2-1}  \\
&& [{\hat J}_{3},{\hat J}_{\pm}]= \pm {\hat J}_{\pm},
   \label{eq:deformed-sl2-2}
\end{eqnarray}
where the hat indices help us to distinguish these modified structures 
with respect to the algebra $sl(2)$. In fact, let us define a new basis of 
the algebra subtended by $J_{\pm}$ and $J_3$ as follows:  
\begin{eqnarray}
 && {\hat J}_{+}=J_{+}\;f^{+}({\cal C}, J_{3}),\qquad
  {\hat J}_{-}=f^{-}({\cal C}, J_{3})\;J_{-},
 \label{eq:fun-real}
\end{eqnarray}
and
\begin{eqnarray}
 && {\hat J}_{3}=J_{3}, \qquad\qquad 
 \label{eq:fun-real-J3}
\end{eqnarray}
so that we evidently ensure the relations (\ref{eq:deformed-sl2-1}) and 
(\ref{eq:deformed-sl2-2}) for 
arbitrary functions $f^{+}$ and $f^{-}$ in terms of the commuting operators 
${\cal C}$ and $J_{3}$ if we require that, on the state $|j,m\rangle$, 
we have
\begin{eqnarray}
 && (j+m)(j-m+1)f^{+}(j,m-1) f^{-}(j,m-1)- \nonumber \\
 && (j-m)(j+m+1)f^{+}(j,m) f^{-}(j,m) 
 = \sum_{p=0}^{N}\beta_{p}(2m)^{2p+1}.
 \label{eq:diff-equ}
\end{eqnarray}
If $f^{\pm}$ are real functions of ${\cal C}$ and $J_{3}$, then hermiticity 
implies $f^{+}=f^{-}$.

Let us point out that our choice (\ref{eq:fun-real}) is such that 
the ladder generators can be seen as hermitian conjugate ones and that 
eq. (\ref{eq:fun-real-J3}) leaves unchanged the diagonal operator 
$J_{3}$. Relatively fastidious calculations starting with eq. 
(\ref{eq:diff-equ}) lead to the result
\begin{eqnarray}
  && (j-m)(j+m+1)f^{+}(j,m)f^{-}(j,m) =\sum_{p=0}^{N}\beta_{p}2^{2p+1}
   \biggl(\sum_{r=1}^{j}r^{2p+1}-\sum_{r=1}^{m}r^{2p+1}\biggr)\nonumber \\
   && =\sum_{p=0}^{N}\beta_{p} 2^{2p+1} \biggl( {1\over 2p+2}j^{2p+2}+
  {1\over 2}j^{2p+1}+{1\over 2}\biggl({2p+1 \atop 1}\biggr) B_{1}j^{2p} 
  -{1\over 4}\biggl({2p+1 \atop 3}\biggr)B_{2}j^{2p-2} \nonumber\\
  && +{1\over 6}\biggl({2p+1 \atop 5}\biggr)B_{3}j^{2p-4}-\cdots- 
     (-1)^{p}{1\over 2p}\biggl({2p+1 \atop 2p-1}\biggr)B_{p}j^{2}-
{1\over 2p+2}m^{2p+2}- {1\over 2}m^{2p+1} \nonumber\\
  && -{1\over 2}\biggl({2p+1 \atop 1}\biggr) B_{1}m^{2p} 
  +{1\over 4}\biggl({2p+1 \atop 3}\biggr)B_{2}m^{2p-2} 
  -{1\over 6}\biggl({2p+1 \atop 5}\biggr)B_{3}m^{2p-4}+\cdots+ \nonumber\\
  && (-1)^{p}{1\over 2p}\biggl({2p+1 \atop 2p-1}\biggr)B_{p}m^{2} )\biggr) 
\end{eqnarray}
where $B_{1}={1 \over 6}$, $B_{2}={1 \over 30}$, $B_{3}= {1 \over 42}$,
$\cdots$, are Bernoulli numbers \cite{Jolley} appearing in this particular 
summation of series \cite{Jolley}. By dividing both sides of the above 
equality by 
$(j-m)(j+m+1)$, the final result can be put in the form
\begin{eqnarray}
f^{+}(j,m)f^{-}(j,m)=\beta_{0}+\sum_{k=1}^{N}\beta_{k}{2^{2k} \over k+1}
\biggl( \sum_{r=1}^{k} \sum_{s=0}^{r}(j(j+1))^{s}(m(m+1))^{r-s}\epsilon_{r}(k)
\biggr)
\label{eq:f-exp}
\end{eqnarray}
or, in terms of generators,
\begin{eqnarray}
f^{+}({\cal C},J_{3})f^{-}({\cal C},J_{3})=
\beta_{0}+\sum_{k=1}^{N}\beta_{k}{2^{2k} \over k+1}
\biggl(\sum_{r=1}^{k} \epsilon_{r}(k)\sum_{s=0}^{r}{\cal C}^{s}
(J_{3}(J_{3}+1))^{r-s}\biggr).
\label{eq:f-oper}
\end{eqnarray}
In eqs. (\ref{eq:f-exp}) and (\ref{eq:f-oper}), we have introduced specific 
functions of $k$ defined by the following relations
\begin{eqnarray}
\epsilon_{k}(k)=1,
\label{eq:def-epsilon-1}
\end{eqnarray} 
and, for $j=1,\;2,\cdots,\;k-1$,
\begin{eqnarray}
  (-1)^{j+1}\biggl({k+1\over j}\biggr)\biggl({2k+1\atop 2j-1}\biggr)B_{j} &=&
\biggl({k+1\atop 2j}\biggr)+\epsilon_{k-1}(k)\biggl({k\atop 2j-2}\biggr) 
\nonumber \\
   && +\epsilon_{k-2}(k)\biggl({k-1\atop 2j-4}\biggr)+\cdots+\epsilon_{k-j}(k).
\label{eq:def-epsilon-2}
\end{eqnarray} 

In addition, let us also point out that we could rewrite eq. (\ref{eq:f-oper})
on the following form
\begin{eqnarray}
f^{+}({\cal C},J_{3})f^{-}({\cal C},J_{3})=
\sum_{k=1}^{N+1}\alpha_{k}\biggl(
\sum_{n=0}^{k-1} {\cal C}^{k-1-n}
(J_{3}(J_{3}+1))^{n}\biggr),
\end{eqnarray}
leading to simple identifications between the $\alpha$- and 
$\beta$-coefficients. In fact, we have
\begin{eqnarray}
  \alpha_{1}=\beta_{0}, \qquad\qquad \alpha_{l}=\sum_{k=l-1}^{N}\beta_{k}
  {2^{2k} \over k+1} \epsilon_{l-1}(k),\qquad\qquad l=2,\;3,\cdots,\;N+1.
\label{eq:alpha(beta)}
\end{eqnarray}     
With this last set of information, the corresponding representations are 
simpler. Indeed, we get 
\begin{eqnarray}
{\hat J}_{\pm}|j,m\rangle =\biggl( \sum_{k=1}^{N+1}\alpha_{k}
\bigl((j(j+1))^{k}-(m(m+1))^{k}\bigr)\biggr)^{1/2}|j,m\pm 1\rangle,
\label{eq:J-action}
\end{eqnarray}
and the commutation relation (\ref{eq:deformed-sl2-1}) becomes
\begin{eqnarray}  
  [{\hat J}_{+}, {\hat J}_{-}] = 2 \sum_{n=1}^{N+1}\alpha_{n}
\sum_{r=0}^{R_{n}}\biggl({n \atop 2r+1}\biggr){\hat J}_{3}^{2n-2r-1},
\label{eq:commutator-1}
\end{eqnarray}
where $R_{n}={1\over 2}(n-2)$ for even $n$ and $R_{n}={1\over 2}(n-1)$ for 
odd $n$. We thus relate the $\alpha$- and $\beta$-coefficients in the 
other way (with respect to eqs. (\ref{eq:alpha(beta)})) by
\begin{eqnarray}
\beta_{p}=2^{-2p}\sum_{k=p+1}^{2p+1}\alpha_{k}\biggl({k \atop 2k-2p-1}\biggr),
\qquad\qquad  p=0,\;1,\cdots,\;N.
\label{eq:beta(alpha)}
\end{eqnarray}

Up to these choices, we have obtained at this stage some new information 
on irreducible representations of the nonlinear $sl(2)$ algebra for $N$ 
arbitrary. We have to add more specific arguments in order to get all the 
representations as it will appear in the following.

Now, let us give the explicit expressions of the deformed generators 
${\hat J}_{\pm}$. According to eqs. (\ref{eq:fun-real}) and
(\ref{eq:fun-real-J3}), we have
\begin{eqnarray}
  && {\hat J}_{+}=J_{+}\biggl( \sum_{k=1}^{N+1}\sum_{r=0}^{k-1}
\alpha_{k} {\cal C}^{k-1-r}(J_{3}(J_{3}+1))^{r}\biggr)^{1/2},
\end{eqnarray}
and
\begin{eqnarray}
  && {\hat J}_{-}=\biggl( \sum_{k=1}^{N+1}\sum_{r=0}^{k-1}
\alpha_{k} {\cal C}^{k-1-r}(J_{3}(J_{3}+1))^{r}\biggr)^{1/2}J_{-},
\end{eqnarray}
or
\begin{eqnarray}
  && {\hat J}_{+}=J_{+}\biggl( \sum_{k=1}^{N+1}
\alpha_{k} {{\cal C}^{k}-(J_{3}(J_{3}+1))^{k} \over 
{\cal C}-J_{3}(J_{3}+1)} \biggr)^{1/2},
\end{eqnarray}
and
\begin{eqnarray}
  && {\hat J}_{-}=\biggl( \sum_{k=1}^{N+1}
\alpha_{k} {{\cal C}^{k}-(J_{3}(J_{3}+1))^{k} \over 
{\cal C}-J_{3}(J_{3}+1)} \biggr)^{1/2}J_{-}.
\end{eqnarray}
Moreover, if we define
\begin{eqnarray}
\phi(x)=\sum_{k=1}^{N+1}\alpha_{k}\;x^{k},
\label{eq:phi}
\end{eqnarray}
these generators become
\begin{eqnarray}
  && {\hat J}_{+}=J_{+}\biggl( 
{\phi({\cal C})-\phi(J_{3}(J_{3}+1)) \over 
{\cal C}-J_{3}(J_{3}+1)} \biggr)^{1/2}, \label{eq:deformed(undeformed)-1} \\
  && \nonumber \\
  && {\hat J}_{-}=\biggl( {\phi({\cal C})-\phi(J_{3}(J_{3}+1)) \over 
     {\cal C}-J_{3}(J_{3}+1)} \biggr)^{1/2} J_{-},
      \label{eq:deformed(undeformed)-2}
\end{eqnarray}
and the corresponding Casimir operator is
\begin{eqnarray}
  && {\hat {\cal C}}={1 \over 2}\biggl({\hat J}_{+}{\hat J}_{-}+
     {\hat J}_{-}{\hat J}_{+}+\phi({\hat J}_{3}({\hat J}_{3}+1))+
     \phi({\hat J}_{3}({\hat J}_{3}-1))\biggr)=\phi({\cal C}).
  \label{eq:deformed-Casimir}
\end{eqnarray}

\noindent {\bf Remark:} We can can also write (\ref{eq:deformed(undeformed)-1}) and 
(\ref{eq:deformed(undeformed)-2}) as (see the ${\cal U}_{q}(sl(2))$ case)
\begin{eqnarray}
  && {\hat J}_{+}=J_{+}\Biggl( {\biggl(\psi\biggl(\sqrt{{\cal C}+
\biggl({1\over 2}\biggr)^{2}}\biggr)\biggr)^{2}-
\biggl(\psi(J_{3}+{1\over 2})\biggr)^{2} \over 
\biggl(\sqrt{{\cal C}+\biggl({1\over 2}\biggr)^{2}}\biggr)^{2}-
\biggl(J_{3}+{1\over 2}\biggr)^{2}} \Biggr)^{1/2}, 
\label{eq:deformed(undeformed)-1-2} \\
  && \nonumber \\
  && {\hat J}_{-}=\Biggl( {\biggl(\psi\biggl(\sqrt{{\cal C}+
\biggl({1\over 2}\biggr)^{2}}\biggr)\biggr)^{2}-
\biggl(\psi(J_{3}+{1\over 2})\biggr)^{2} \over 
\biggl(\sqrt{{\cal C}+\biggl({1\over 2}\biggr)^{2}}\biggr)^{2}-
\biggl(J_{3}+{1\over 2}\biggr)^{2}} \Biggr)^{1/2}J_{-},
      \label{eq:deformed(undeformed)-2-2}
\end{eqnarray}
where, 
\begin{eqnarray}
\phi(x)=\psi^{2}\biggl(\sqrt{x+{1\over 4}}\biggr)-
\psi^{2}\biggl({1\over 2}\biggr),\qquad\hbox{if}\qquad \phi(0)=0.
\end{eqnarray}
The relation between the deformed Casimir ${\hat {\cal C}}$ and 
${\cal C}$ is given by
\begin{eqnarray}
\sqrt{{\hat {\cal C}}+\biggl(\psi\biggl({1\over 2}\biggr)\biggr)^{2}}=
\psi\biggl(\sqrt{{\cal C}+\biggl({1\over 2}\biggr)^{2}}\biggr).
\end{eqnarray}

Now, if $\phi$ is bijective, we evidently have
\begin{eqnarray}
{\cal C}=\phi^{-1}({\hat {\cal C}}),
\label{eq:bijective-condition}
\end{eqnarray}
and
\begin{eqnarray}
  &&  J_{+}={\hat J}_{+}\biggl( 
{\phi^{-1}({\hat {\cal C}})-{\hat J}_{3}({\hat J}_{3}+1) \over 
{\hat {\cal C}}-\phi({\hat J}_{3}({\hat J}_{3}+1))} \biggr)^{1/2},
\label{eq:undeformed(deformed)-1}  \\
  && \nonumber \\
  &&  J_{-}=\biggl( 
{\phi^{-1}({\hat {\cal C}})-{\hat J}_{3}({\hat J}_{3}+1) \over 
{\hat {\cal C}}-\phi({\hat J}_{3}({\hat J}_{3}+1))} \biggr)^{1/2}{\hat J}_{-}.
\label{eq:undeformed(deformed)-2}
\end{eqnarray}
From this point of view bijective $\phi$ ' s are of particular 
interest. A similar discussion is valid for the function $\psi$.

\sect{A Generalization\label{sect:generalisation}}

For (\ref{eq:deformed(undeformed)-1}) and (\ref{eq:deformed(undeformed)-2}), 
the starting point is $sl(2)$. One can take ${\cal U}_{q}(sl(2))$ (which 
is itself a nonlinear generalization of $sl(2)$) as the starting point and 
generalize that again by postulating
\begin{eqnarray}
{\hat J}_{\pm}|j,m\rangle =\biggl( \sum_{k=1}^{N+1}\alpha_{k}
\bigl(([j][j+1])^{k}-([m][m+1])^{k}\bigr)\biggr)^{1/2}|j,m\pm 1\rangle,
\label{eq:q-J-action}
\end{eqnarray}
i.e.
\begin{eqnarray}
  && {\hat J}_{+}=J_{+}\biggl({\phi({\cal C})-\phi([J_{3}][J_{3}+1]) \over 
	{\cal C}-[J_{3}][J_{3}+1]} \biggr)^{1/2}, 
	\label{eq:deformed(undeformed)-1-q} \\
  && \nonumber \\
  && {\hat J}_{-}=\biggl( {\phi({\cal C})-\phi([J_{3}][J_{3}+1]) \over 
     {\cal C}-[J_{3}][J_{3}+1]} \biggr)^{1/2} J_{-},
      \label{eq:deformed(undeformed)-2-q}
\end{eqnarray}  
where
\begin{eqnarray}
{\cal C}={1\over 2}(J_{+}J_{-}+J_{-}J_{+})+[J_{3}]^{2}.
\end{eqnarray}  
For example, if we choose 
\begin{eqnarray}
\phi(x)=x+{\beta \over [2]}x^{2},
\end{eqnarray} 
we obtain the following commutation relation
\begin{eqnarray}
[{\hat J}_{+}, \;{\hat J}_{-}]=[2J_{3}](1+\beta [J_{3}]^{2}).
\end{eqnarray} 
For another choice, we can obtain
\begin{eqnarray}
[{\hat J}_{+}, \;{\hat J}_{-}]=[\;[2J_{3}]_{1}]_{2},
\label{eq:generalisation-1}
\end{eqnarray}  
where,
\begin{eqnarray}
[x]_{i}={q_{i}^{x}-q_{i}^{-x} \over q_{i}-q_{i}^{-1}}, \qquad q_{i}\in \CC.
\end{eqnarray}  
One can ultimately even envisage a hierarchy of $q$-brackets 
generalizing the r.h.s. of (\ref{eq:generalisation-1}). 

When one generalizes (\ref{eq:J-action}) as in (\ref{eq:q-J-action}), 
(${\hat J}_{\pm}$, $q^{\pm {\hat J}_{3}}$) being expressed in terms of 
($J_{\pm}$, $q^{\pm  J_{3}}$) of ${\cal U}_{q}(sl(2))$, one can implement 
the standard Hopf structure of the latter (rather than starting from 
that of $sl(2)$) to construct ${\hat J}_{\pm}$ for product representations.
Evidently the formalism of this section contains the results of the 
preceding one as limiting cases ($q \longrightarrow 1$).

\sect{The $sl(2)$ and ${\cal U}_{q}(sl(2))$-Contexts\label{sect:sl2q-case}}

The nonlinear $sl(2)$-algebras given by eqs. (\ref{eq:deformed-sl2-1}) and
(\ref{eq:deformed-sl2-2}) evidently contain the 
expected linear $sl(2)$-one as well as its $q$-deformation 
${\cal U}_{q}(sl(2))$. The first one corresponds to $N=0$, so that eqs. 
(\ref{eq:deformed-sl2-1}) with $\beta_{0}=1$ and (\ref{eq:def-sl2-1}) become 
identical while the second one is readily obtained by taking the limit 
$N\rightarrow  \infty$ with the coefficients
\begin{eqnarray}  
  \beta_{p}&=&{2 \over q-q^{-1}}{(\log q)^{2p+1} \over (2p+1)!}, 
\qquad\qquad p=0,\;1,\cdots \label{eq:q-beta-1} \\ 
  &=&{1 \over \sinh \delta}{(\delta)^{2p+1} \over (2p+1)!},
\qquad\qquad q=\exp\;\delta. 
\label{eq:q-beta-2}
\end{eqnarray}
If, in the linear case, we evidently have
\begin{eqnarray} 
  && f^{+}(j,m)f^{-}(j,m)=1,\qquad\qquad\qquad\qquad f^{+}=f^{-}=1,
\end{eqnarray}
ensuring that
\begin{eqnarray}  
{\hat J}_{\pm}=J_{\pm},\qquad\qquad\qquad\qquad  {\hat J}_{3}=J_{3},
\end{eqnarray}
we point out in the $q$-deformation that \cite{Zachos}
\begin{eqnarray}  
&& f^{+}(j,m)f^{-}(j,m)={ [j-m][j+m+1] \over (j-m)(j+m+1)},
\label{eq:deformed-map}
\end{eqnarray}
where as usual
\begin{eqnarray}  
[x]={q^{x}-q^{-x} \over q-q^{-1}}.
\label{eq:q-number}
\end{eqnarray}
By developing the right-hand side of eq. (\ref{eq:deformed-map}), it is not 
difficult to show that it coincides with our expression (\ref{eq:f-exp}) 
(for example) with the coefficients (\ref{eq:q-beta-1}). This corresponds 
to the equality
\begin{eqnarray}  
  && {1 \over (j-m)(j+m+1)}{\cosh \bigl(\delta (2j+1)\bigr)- 
\cosh \bigl(\delta(2m+1)\bigr) \over
     2 \sinh^{2}\delta} \nonumber \\
  && ={\delta \over \sin \delta}+\sum_{k=1}^{\infty}{2^{2k+1}\delta^{2k+1}
     \over (2k+2)!\;\sinh \delta}\sum_{r=1}^{k}\sum_{s=0}^{r}(j(j+1))^{s}
     (m(m+1))^{r-s}\epsilon_{r}(k) 
\end{eqnarray} 
where the corresponding functions $\epsilon_{r}(k)$ are given by 
(\ref{eq:def-epsilon-1}) and (\ref{eq:def-epsilon-2}).

Let us also mention that the quantum algebra ${\cal U}_{q}(sl(2))$
corresponds to the choice of the following bijective function 
introduced by eq. (\ref{eq:phi}):
\begin{eqnarray} 
\phi (J_{3}(J_{3}+1)) = [J_{3}][J_{3}+1],\qquad\qquad
\phi({\cal C})=[\sqrt{{\cal C}+{1 \over 4}}-{1\over 2}]
[\sqrt{{\cal C}+{1 \over 4}}+{1\over 2}],
\end{eqnarray} 
with the bracket (\ref{eq:q-number}) so that the generators 
(\ref{eq:deformed(undeformed)-1}) and (\ref{eq:deformed(undeformed)-2})  
become in this context
\begin{eqnarray} 
 &&  {\hat J}_{+}=J_{+}\Biggl({\bigl[\sqrt{{\cal C}+
\biggl({1 \over 2}\biggr)^{2}}\bigr]^{2}-
\bigl[J_{3}+{1\over 2}\bigr]^{2} \over 
\biggl(\sqrt{{\cal C}+\biggl({1 \over 2}\biggr)^{2}}\biggr)^{2}-
\biggl(J_{3}+{1\over 2}\biggr)^{2}}\Biggr)^{1/2},
\label{eq:q-generator-1}  \\
 && {\hat J}_{-}= \Biggl({\bigl[\sqrt{{\cal C}+
\biggl({1 \over 2}\biggr)^{2}}\bigr]^{2}-
\bigl[J_{3}+{1\over 2}\bigr]^{2} \over 
\biggl(\sqrt{{\cal C}+\biggl({1 \over 2}\biggr)^{2}}\biggr)^{2}-
\biggl(J_{3}+{1\over 2}\biggr)^{2}}\Biggr)^{1/2}J_{-}.
\label{eq:q-generator-2} 
\end{eqnarray} 
The corresponding Casimir operator is then given by
\begin{eqnarray} 
{\hat {\cal C}}=\bigl[\sqrt{{\cal C}+\biggl({1 \over 2}\biggr)^{2}}-
{1\over 2}\bigr]\bigl[\sqrt{{\cal C}+\biggl({1 \over 2}\biggr)^{2}}+
{1\over 2}\bigr],
\label{eq:q-Casimir}
\end{eqnarray}
i.e.
\begin{eqnarray} 
\sqrt{{\hat {\cal C}}+\biggl[{1\over 2}\biggr]^{2}}
=\biggl[\sqrt{{\cal C}+\biggl({1 \over 2}\biggr)^{2}}\biggr],
\end{eqnarray}
and, consequently,
\begin{eqnarray}
\sqrt{{\cal C}+\biggl({1\over 2}\biggr)^{2}}=
 {1\over \delta}\;\arcsinh\biggl(\sqrt{{\hat {\cal C}}+
\biggl[{1 \over 2}\biggr]^{2}}\sinh\delta \biggr) .
\label{eq:q-Casimir-bijective}
\end{eqnarray} 
These relations finally lead to
\begin{eqnarray} 
 && J_{+}={\hat J}_{+}\Biggl({\biggl( {1\over \delta}\;\arcsinh(\sqrt{{\hat 
    {\cal C}}+\biggl[{1 \over 2}\biggr]^{2}}\;\sinh\delta ) \biggr)^{2}-
   ({\hat J}_{3}+{1\over 2})^{2}\over \biggl(\sqrt{{\hat {\cal C}}+
    \biggl[{1 \over 2}\biggr]^{2}} \biggr)^{2}-[{\hat J}_{3}+
   {1\over 2}]^{2}}\Biggr)^{1/2}, 
     \label{eq:q-generator-1-bijective} \\
  && J_{-}=(J_{+})^{{\cal +}}=\Biggl({\biggl( {1\over \delta}\;\arcsinh(\sqrt{{\hat 
   {\cal C}}+\biggl[{1 \over 2}\biggr]^{2}}\;\sinh\delta) \biggr)^{2}-
  ({\hat J}_{3}+{1\over 2})^{2}\over \biggl(\sqrt{{\hat {\cal C}}+ 
     \biggl[{1 \over 2}\biggr]^{2}} \biggr)^{2}-
    [{\hat J}_{3}+{1 \over 2}]^{2}}\Biggr)^{1/2}{\hat J}_{-}, 
   \label{eq:q-generator-2-bijective}
\end{eqnarray}
ensuring that we have the expected commutation relation
\begin{eqnarray} 
  [{\hat J}_{+}, {\hat J}_{-}]= [2J_{3}].
\label{eq:q-commutator}
\end{eqnarray}

It has been shown by Curtright et al \cite{Zachos} that, from the 
well known $sl(2)$-co-commutative coproduct and eqs. 
(\ref{eq:q-generator-1})-(\ref{eq:q-generator-2}), it is possible to 
characterize ${\cal U}_{q}(sl(2))$ by a co-commutative coproduct.
Moreover,  by the inverse map 
(\ref{eq:q-generator-1-bijective})-(\ref{eq:q-generator-2-bijective}) and
the nonco-commutative coproduct of ${\cal U}_{q}(sl(2))$ noted by 
$\bigtriangleup_{q}$, we can also characterize the linear $sl(2)$ by a 
nonco-commutative one.   

We do not go further into these directions due to our specific interest 
in {\sl finite} values of $N\not = 0$ and more particularly in $N=1$, the
first nontrivial value which has a direct connection with already 
studied physical contexts \cite{Zhedanov}.

\sect{Hopf Structure of Nonlinear Algebras
\label{sect:Hopf-structure}}

In the following, we start by enlarging the term of enveloping algebra of 
$sl(2)$ to include square roots. Then, exploiting the well known fact 
\cite{DJ} that the undeformed generators $J_{\pm}$ and $J_{3}$ admit a 
Hopf structure with the well-known coproduct, counit and antipode given 
for example in the {\sl co-commutative case} respectively by
\begin{eqnarray}
(i)\qquad\qquad  && \bigtriangleup (J_{\pm})=J_{\pm}\otimes 1 + 1 \otimes J_{\pm},\qquad\qquad \qquad\qquad 
	\label{eq:sl2-Hopf-1} \\
  && \nonumber \\ 
  && \bigtriangleup (J_{3})=J_{3}\otimes 1 + 1 \otimes J_{3},
	\label{eq:sl2-Hopf-2} 
\end{eqnarray}
leading to
\begin{eqnarray}
 && \bigtriangleup ({\cal C})={\cal C}\otimes 1 + 1 \otimes {\cal C}+
J_{+}\otimes J_{-} + J_{-}\otimes J_{+}+2 J_{3}\otimes J_{3},
	\label{eq:sl2-Hopf-3}  \\
 && \nonumber \\
(ii)\qquad  && \varepsilon (J_{\pm})=\varepsilon (J_{3})=\varepsilon ({\cal C})=0,
\qquad\qquad 
	\label{eq:sl2-Hopf-4}  \\
 && \nonumber \\
(iii)\qquad  && S(J_{\pm})=-J_{\pm},\qquad\qquad S(J_{3})=-J_{3}, \qquad\qquad S(C)=C, \qquad\qquad 
	\label{eq:sl2-Hopf-5}, 
\end{eqnarray}
we can deduce that our deformed generators ${\hat J}_{\pm}$ and 
${\hat J}_{3}$ also satisfy the Hopf axioms, i.e. \cite{DJ}:
\begin{eqnarray}
  && (\hbox{id}\otimes \bigtriangleup)\bigtriangleup =
      (\bigtriangleup \otimes \hbox{id})\bigtriangleup, 
	\label{eq:sl2-Hopf-6} \\
  && \nonumber \\
  && m(\hbox{id}\otimes S)\bigtriangleup= 
      m(S \otimes \hbox{id})\bigtriangleup= i\circ\varepsilon,
	\label{eq:sl2-Hopf-7}  \\
  && \nonumber \\
  && (\hbox{id}\otimes \varepsilon )\bigtriangleup =
      (\varepsilon \otimes \hbox{id})\bigtriangleup=\hbox{id}.
	\label{eq:sl2-Hopf-8} 
\end{eqnarray}

Now, the coproduct of our deformed generators is given by
\begin{eqnarray}
  && \bigtriangleup ({\hat J}_{+})= f^{+}\biggl( \bigtriangleup ({\cal C}),
\bigtriangleup (J_{3})\biggr)\bigtriangleup ( J_{+}), 
\label{eq:deformed-coproduct-1}\\
  && \nonumber \\
  && \bigtriangleup ({\hat J}_{-})=\bigtriangleup ( J_{-})
  f^{-}\biggl( \bigtriangleup ({\cal C}),\bigtriangleup (J_{3})\biggr),
 \label{eq:deformed-coproduct-2}
\end{eqnarray}
it is not difficult to test that this coproduct is co-commutative,
the same way of reasoning applying to the counit and antipode. 

Let us remark that the r.h.s. of (\ref{eq:deformed-coproduct-1}) and 
(\ref{eq:deformed-coproduct-2}) is just an expansion of the generators 
$J_{\pm}$ and $J_{3}$. If the function $\phi$ ($\psi$) is bijective
eqs. (\ref{eq:undeformed(deformed)-1})-(\ref{eq:undeformed(deformed)-2}), 
the Hopf structure {\sl can be written} using {\sl only} the {\sl deformed} 
generators ${\hat J}_{\pm}$ and ${\hat J}_{3}$. If we take 
$\bigtriangleup_{q}$ given by Curtright et al \cite{Zachos} we can endow 
the nonlinear algebra ($\phi$ is bijective) by a nonco-commutative 
coproduct.

\sect{The Cubic $sl(2)$-Algebra\label{sect:cubic-case}}

Let us now consider the $N=1$-context leading, in eq. 
(\ref{eq:deformed-sl2-1}), at most to 
the {\sl cubic} power in the diagonal generator. This corresponds in 
particular to the nonlinear Higgs algebra \cite{Higgs}, a symmetry one
for the harmonic oscillator and the Kepler problems in a two-dimensional 
curved space. From eq. (\ref{eq:f-exp}), we immediately get
\begin{eqnarray}  
  f^{+}(j,m)f^{-}(j,m)=\beta_{0}+2\beta_{1}\biggl(j(j+1)+m(m+1)\biggr),  
  \label{eq:cubic-exp}
\end{eqnarray}
leading to the Higgs algebra when $\beta_{0}=1$, $\beta_{1}=\beta$. Then, 
we have the commutation relations
\begin{eqnarray}
  && [{\hat J}_{+}, {\hat J}_{-}]= 2{\hat J}_{3}+8 \beta {\hat J}_{3}^{3}, 
     \label{eq:cubic-sl2-1} \\
  && \nonumber \\
  && [ {\hat J}_{3}, {\hat J}_{\pm}]= \pm {\hat J}_{\pm}.
  \label{eq:cubic-sl2-2}
\end{eqnarray}
By requiring that the ladder operators are hermitian conjugate to each other,
we have to fix
\begin{eqnarray}
  f^{+}(j,m)=\biggl(1+2\beta \biggl(j(j+1)+m(m+1)\biggr)\biggr)^{1/2},
  \label{eq:f-cubic}
\end{eqnarray}
so that the unitary irreducible representations of the Higgs algebra are 
given by
\begin{eqnarray}
  && {\hat J}_{\pm}|j,m\rangle = \biggl(
(j\mp m)(j\pm m+1)\biggl(1+2\beta (j(j+1)+m(m\pm 1))\biggr)\biggr)^{1/2}
|j,m\pm 1 \rangle,\qquad \label{eq:cubic-rep-1}\\
  && \nonumber \\
  && {\hat J}_{3}|j,m\rangle = m |j,m\rangle,
  \label{eq:cubic-rep-2}
\end{eqnarray}
where the parameter $\beta$ is constrained by ensuring 
\begin{eqnarray}
1+2\beta \biggl(j(j+1)+m(m\pm 1)\biggr) \geq 0,
\end{eqnarray}
or 
\begin{eqnarray}
\beta \geq -{1 \over 4 j^{2}}, \qquad\qquad \forall \;j\;\;(j\not = 0).
\label{eq:beta-condition} 
\end{eqnarray}

Such unitary irreducible representations (\ref{eq:cubic-rep-1}) and
(\ref{eq:cubic-rep-2}) are associated with explicit forms of the 
$sl_{\beta}(2)$-nonlinear generators expressed in terms of the undeformed 
$sl(2)$-ones. In fact, we formally claim that, according to eqs. 
(\ref{eq:fun-real}) and (\ref{eq:cubic-rep-1}), we have
\begin{eqnarray}
  {\hat J}_{+}=J_{+}\biggl(1+2\beta ({\cal C}+J_{3}(J_{3}+1))\biggr)^{1/2}
={\cal Q}({\cal C}, J_{+}, J_{3}), 
\end{eqnarray}
and
\begin{eqnarray}
  {\hat J}_{-}=\biggl(1+2\beta ({\cal C}+J_{3}(J_{3}+1))\biggr)^{1/2}J_{-}
=({\hat J}_{+})^{+}, 
\end{eqnarray}
while the third one ${\hat J}_{3}$ is unchanged (see eq. 
(\ref{eq:fun-real-J3})). Here we point out that the corresponding 
$\phi$-function (\ref{eq:phi}) is not bijective and the corresponding Hopf
structure cannot be written using our $\beta$-generators.

\hfill

Let us now insist on an interesting property which, at our knowledge, 
seems to be not yet exploited, i.e. on a possible shift of the diagonal 
generator spectrum expressed in terms of a (real scalar) parameter called 
hereafter $\gamma$. So, let us propose to modify the relation 
(\ref{eq:cubic-rep-2}) in the following way
\begin{eqnarray}     
  {\hat J}_{3}|j,m\rangle = (m+\gamma)|j,m\rangle.
 \label{eq:shift}
\end{eqnarray}
If it is evident that, in the usual angular momentum theory, such a shift 
has no physical meaning, it is not trivial to show that, in a 
$q$-deformed one, nothing more happens when $q$ is not a root of 
unity. Indeed, if we ask for the commutation relation
\begin{eqnarray}  
[J_{+}, J_{-}]=[2J_{3}]
\end{eqnarray}
with the bracket (\ref{eq:q-number}) and if we require
\begin{eqnarray}   
  && J_{+}|j,m\rangle=\sqrt{f(j,m)}\;|j,m+1\rangle, \label{eq:xx-1}\\
  && \nonumber \\
  && J_{+}|j,m\rangle=\sqrt{f(j,m-1)}\;|j,m-1\rangle,
  \label{eq:xx-2}
\end{eqnarray}
when
\begin{eqnarray}   
  J_{3}|j,m\rangle=(m+\gamma)|j,m\rangle, \label{eq:xx-3}
\end{eqnarray}
it is possible to show that
\begin{eqnarray} 
f(j,m)={1 \over (q-q^{-1})^{2}}\biggl(q^{-2j+2\gamma -1}+q^{2j-2\gamma +1}-
q^{2m+2\gamma +1}-q^{-2m-2\gamma -1}\biggr). \label{eq:zz}
\end{eqnarray}
Then, due to the fact that, from eq. (\ref{eq:xx-1}), we have
\begin{eqnarray} 
f(j,j)=0,
\end{eqnarray}
we get from eq. (\ref{eq:zz})
\begin{eqnarray} 
f(j,j)=[-2\gamma][2j+1]=0,
\end{eqnarray}
asking for the annulation of the parameter $\gamma$. We thus 
conclude that the shift (\ref{eq:xx-3}) does not permit to characterize 
new representations of ${\cal U}_{q}(sl(2))$.
 
The study of the Higgs algebra in that direction is richer and nonzero 
values of $\gamma$ can be exploited in order to select new unitary 
irreducible representations of this cubic $sl(2)$-algebra. In order to 
establish such a result, let us consider eq. (\ref{eq:shift}) inside the Higgs 
context characterized by the commutation relations (\ref{eq:cubic-sl2-1}) and
(\ref{eq:cubic-sl2-2}). The action of the ladder operators ${\hat J}_{\pm}$ 
on the basis leads to $\gamma$-dependent $f^{\pm}$-functions. In fact, 
in correspondence with eqs. (\ref{eq:cubic-rep-1}), we get here
\begin{eqnarray}       
{\hat J}_{+}|j,m\rangle &=&\biggl((j-m)(j+m+1+2\gamma)\biggl(1+
2\beta(j(j+1) \nonumber \\
&&+m(m+1)+2\gamma(j+m+1+\gamma))\biggr)\biggr)^{1/2}|j,m+1\rangle,
\label{Higgs-rep-shift-1}
\end{eqnarray}
and
\begin{eqnarray}       
{\hat J}_{-}|j,m\rangle &=& \biggl((j-m+1)(j+m+2\gamma)\biggl(1+
2\beta(j(j+1)   \nonumber\\
&& +m(m-1)+2\gamma(j+m+\gamma))\biggr)\biggr)^{1/2}|j,m-1\rangle.
\label{Higgs-rep-shift-2}
\end{eqnarray}

By exploiting the property that
\begin{eqnarray}  
{\hat J}_{-}|j,-j\rangle =0,
\end{eqnarray} 
we get the constraint
\begin{eqnarray}  
2\gamma(2j+1)\biggl(1+4\beta(j(j+1)+\gamma^2)\biggl)=0, 
\end{eqnarray} 
showing that, besides our preceding context ($\gamma=0$), there are other 
possibilities related to nonzero $\gamma$-values issued from the equation
\begin{eqnarray}  
  \gamma^{2}={1\over 4\beta^{2}}\biggl(-\beta-4\beta^{2} \;j(j+1)\biggr).
\end{eqnarray} 
A simple discussion of its roots leads to the {\sl two} families of 
new representations characterized respectively by
\begin{eqnarray}  
  \gamma={1\over 2\beta}\biggl(-\beta -4\beta^{2} \;j(j+1)\biggr)^{1/2},
\end{eqnarray}  
or
\begin{eqnarray}  
  \gamma=-{1\over 2\beta}\biggl(-\beta -4\beta^{2} \;j(j+1)\biggr)^{1/2},
\end{eqnarray}  
both values being constrained by the deformation parameter $\beta$ 
such that
\begin{eqnarray}  
-{1\over 4j(j+1)}< \beta \leq -{1\over 4j(j+1)+1}. 
\label{eq:ff}
\end{eqnarray} 

Let us insist on the fact that these representations are typical of the 
deformation characterizing the Higgs algebra: they do not exist when 
$\beta=0$. Moreover, such a method suggests its application to other 
nonlinear $sl(2)$ algebras and we want to look at its impact here on 
an interesting quadratic one \cite{Schoutens} in the following section.

Just as the simplest example, let us fix $j={1\over 2}$ (corresponding to
the fundamental representation in the conventional $sl(2)$-case). We 
evidently conclude that, if our $\beta$-parameter is constrained 
(according to eq. (\ref{eq:ff})) by 
\begin{eqnarray}  
 -{1\over 3}< \beta \leq -{1\over 4}, 
\end{eqnarray}     
we get {\sl three} families of representations corresponding to
\begin{eqnarray}  
  \gamma=\pm {1\over 2\beta}\biggl(-\beta -3\beta^{2}\biggr)^{1/2}
\qquad\hbox{and}\qquad \gamma =0.
\end{eqnarray}  
According to eq. (\ref{eq:beta-condition}) when $\gamma=0$, we have 
$\beta \geq -1$ and we point out that, if $\beta > {1\over 4}$ or if 
$-1 \leq \beta \leq -{1\over 3}$, we get only one family while, evidently, 
if $\beta < -1$, none representation is admissible.

As a last remark, let us notice that the modification effectively 
introduced in eq. (\ref{eq:shift}) through the $\gamma$-parameter does not 
affect our conclusions on the Hopf structure of the Higgs algebra.  

\sect{The Quadratic $sl(2)$-Algebra\label{sect:quadratic-case}}

Another nonlinear $sl(2)$-algebra is the {\sl quadratic} one \cite{Schoutens}
characterized by the following commutation relations depending on the (real
scalar) parameter $\alpha$
\begin{eqnarray} 
  && [J_{+}^{(\alpha)}, J_{-}^{(\alpha)}]=2J_{3}^{\alpha}+4 
     \alpha (J_{3}^{(\alpha)})^{2}, \label{eq:quadratic-sl2-1}\\
  && \nonumber \\
  && [J_{3}^{(\alpha)}, J_{\pm}^{(\alpha)}]= \pm J_{\pm}^{(\alpha)}.
     \label{eq:quadratic-sl2-2}
\end{eqnarray}
It has already been exploited \cite{Schoutens} in connection with Yang-Mills
type gauge theories and with fundamental quantum mechanical problems 
\cite{Rocek} \cite{Quesne}. In particular, its representation theory has 
already been investigated \cite{Rocek} for the lowest eigenvalues 
of the Casimir operator.

Let us here come back on this representation theory when combined 
with the demand corresponding to eq. (\ref{eq:shift}) of the preceding 
section, i.e.
\begin{eqnarray}  
  J_{3}^{(\alpha)}|j,m\rangle = (m+\gamma)|j,m\rangle.
  \label{eq:quadratic-shift}
\end{eqnarray}
Here the ladder operators $J_{\pm}^{(\alpha)}$ also act on the basis 
and determine $\alpha$-dependent $f^{\pm}$-functions that can be calculated.
They are given in the following relations
\begin{eqnarray}  
 J_{+}^{(\alpha)}|j,m\rangle &=& \biggl((j-m)\biggl(j+m+1+2\gamma+\alpha\bigl(
{4\over 3}j^{2}+ {4\over 3}j m \nonumber \\
  && +{4\over 3}m^{2}+4 \gamma j +4 \gamma m+2j+2m +4\gamma^2 +4 \gamma+
{2\over 3}\bigr)\biggr)\biggr)^{1/2}|j,m+1\rangle,
\label{eq:quadratic-rep-1}
\end{eqnarray}
and
\begin{eqnarray}  
 J_{-}^{(\alpha)}|j,m\rangle &=& \biggl((j-m+1)\biggl(j+m+2\gamma+\alpha\bigl(
{4\over 3}j^{2}+ {4\over 3}j m \nonumber \\
  && +{4\over 3}m^{2}+4 \gamma j +4 \gamma m+{2 \over 3}j-{2 \over 3}m 
+4\gamma^2 \bigr)\biggr)\biggr)^{1/2}|j,m-1\rangle.
\label{eq:quadratic-rep-2}
\end{eqnarray}

Once again, the condition
\begin{eqnarray}  
 J_{-}^{(\alpha)}\;|j,-j\rangle=0, 
\end{eqnarray}
leads to the constraint
\begin{eqnarray}  
\gamma = {1\over 4 \alpha}\Biggl(-1+\sqrt{1-
{16 \over 3}j(j+1)\;\alpha^{2}}\Biggr),
\end{eqnarray}
when
\begin{eqnarray} 
\alpha \leq {3 \over 2(4j+1)}.
\label{eq:alpha-condition}
\end{eqnarray}

Such unitary irreducible representations 
(\ref{eq:quadratic-shift})-(\ref{eq:alpha-condition}) are typical of 
the deformation and are associated with the following forms of 
generators explicitly given in terms of the undeformed $sl(2)$-ones:
\begin{eqnarray} 
&& J_{3}^{(\alpha)}=J_{3}-{1\over 4 \alpha}+{1\over 4 \alpha}
\sqrt{1-{16 \over 3}\;\alpha^{2}\;{\cal C}}, \\
 \label{eq:quadratic-real-1}
&& J_{+}^{(\alpha)}=J_{+}\;\Biggl( {2\over 3}\;\alpha\; (2J_{3}+1)+
\sqrt{1-{16 \over 3}\;\alpha^{2}\;{\cal C}}\Biggr)^{1/2},\\
\label{eq:quadratic-real-2}
&& J_{-}^{(\alpha)}=\Biggl( {2\over 3}\;\alpha \;(2J_{3}+1)+
\sqrt{1-{16 \over 3}\;\alpha^{2}\;{\cal C}}\Biggr)^{1/2}\;J_{-}. 
\label{eq:quadratic-real-3}
\end{eqnarray}
Through the knowledge of the $sl(2)$-coproduct, counit and antipode 
given by eqs. (\ref{eq:sl2-Hopf-1})-(\ref{eq:sl2-Hopf-8}), we can thus 
provide the quadratic algebra
(\ref{eq:quadratic-sl2-1}) and (\ref{eq:quadratic-sl2-2}) with a Hopf 
structure by defining
\begin{eqnarray} 
&& \bigtriangleup(J_{3}^{(\alpha)})=\bigtriangleup(J_{3})-{1\over 4 \alpha}
\;(1\otimes 1)+{1\over 4 \alpha}
\sqrt{1\otimes 1-{16 \over 3}\;\alpha^{2}\;\bigtriangleup({\cal C})},
\label{eq:quadratic-Hopf-1} \\
&& \bigtriangleup(J_{+}^{(\alpha)})=\bigtriangleup(J_{+})
\Biggl( {2\over 3}\;\alpha\;(2\bigtriangleup(J_{3})+1\otimes 1)+
\sqrt{1\otimes 1-{16 \over 3}\;\alpha^{2}
\bigtriangleup({\cal C})}\Biggr)^{1/2},
\label{eq:quadratic-Hopf-2} \\
&& \bigtriangleup(J_{-}^{(\alpha)})=\Biggl( {2\over 3}\;\alpha\; 
(2\bigtriangleup(J_{3})+1\otimes 1 )+
\sqrt{1\otimes 1 -{16 \over 3}\alpha^{2}\bigtriangleup({\cal C})}\Biggr)^{1/2}
\bigtriangleup(J_{-}), 
\label{eq:quadratic-Hopf-3} \\
&& \nonumber \\
&& \varepsilon (J_{3}^{(\alpha)})=\varepsilon (J_{\pm}^{(\alpha)})=0,
\label{eq:quadratic-Hopf-4}  \\
&& \nonumber \\
&& S(J_{3}^{(\alpha)})=-J_{3}-{1\over 4 \alpha}+{1\over 4 \alpha}
\sqrt{1-{16 \over 3}\alpha^{2}{\cal C}},
\label{eq:quadratic-Hopf-5}  \\
&&S(J_{+}^{(\alpha)})=-\Biggl( {2\over 3}\alpha (-2J_{3}+1)+
\sqrt{1-{16 \over 3}\alpha^{2}{\cal C}}\Biggr)^{1/2}J_{+},
\label{eq:quadratic-Hopf-6} \\
&& S(J_{-}^{(\alpha)})=-J_{-}\Biggl( {2\over 3}\alpha (-2J_{3}+1)+
\sqrt{1-{16 \over 3}\alpha^{2}{\cal C}}\Biggr)^{1/2}, 
\label{eq:quadratic-Hopf-7} 
\end{eqnarray} 
as it was the case for the cubic algebra (\ref{eq:cubic-sl2-1}) and
(\ref{eq:cubic-sl2-2}) but with the definitions 
(\ref{eq:sl2-Hopf-1})-(\ref{eq:sl2-Hopf-2}). We note that 
the right-hand side of (\ref{eq:quadratic-Hopf-1})-(\ref{eq:quadratic-Hopf-7})
cannot be written using only the generators $J_{\pm}^{(\alpha)}$ and 
$J_{3}^{(\alpha)}$.

\sect{Conclusions and Comments\label{sect:conclusions}} 

We have just developed the representation theory associated with 
{\sl nonlinear} $sl(2)$-algebras characterized by the structure relations 
(\ref{eq:deformed-sl2-1}) and (\ref{eq:deformed-sl2-2}) containing, in 
particular, the linear $sl(2)$-one as well as its $q$-deformation 
${\cal U}_{q}(sl(2))$. Moreover, we have more specifically visited the 
{\sl cubic} $sl(2)$-algebra in order to get {\sl all} its unitary 
irreducible representations and to show that it is endowed in our 
formalism with a Hopf structure, the corresponding results 
were also presented for the {\sl quadratic} $sl(2)$-algebra. Such a study 
mainly takes advantage of the fact that we can express the generators of 
the nonlinear algebras in terms of the old (undeformed) $sl(2)$-ones and 
that the $sl(2)$-algebra is endowed with a well known Hopf structure. 
These properties allow us to extend our considerations for arbitrary $N$ 
in the odd case (developed in section \ref{sect:representation}) and are 
also valid in principle for the even context after the study of the 
$N=2$-case (developed in section \ref{sect:quadratic-case}).

From the representation point of view, our results generalize to arbitrary 
$j$'s those obtained by Ro${\tilde c}$ek \cite{Rocek}. They also include 
others obtained by Zhedanov \cite{Zhedanov}, Feng Pan \cite{Feng} 
and Bonatsos et al \cite{Bonatsos}.

From the point of view of Hopf structures associated with our developments, 
many connections with recent studies can be pointed out. An interesting
property which has been discussed in section \ref{sect:Hopf-structure} is 
the one concerning the co-commutativity or nonco-commutativity of the 
already known coproducts. We have shown that, in some particular cases, the
nonlinear algebra can be equipped with a consistent Hopf structure
(i.e. the corresponding coproduct being expressed in terms of {\sl deformed}
generators). Moreover, let us mention that there is also a third possibility 
by exploiting our recent proposal for a new deformed structure
${\cal U}_{q}^{\theta}(sl(2))$-algebra using 
a real paragrassmannian variable $\theta$ \cite{Abdess}. 

All these properties have to be carefully examined and we plan to come 
back on these in the future. Let us finally add that our results, in 
particular, confirm  those recently obtained by Quesne and Vansteenkiste    
\cite{Vansteenkiste} showing that if we ask for a deformed coproduct in terms 
of deformed generators, only the already well known ones are possible. We 
have obtained new ones due to the fact that we have expressed the 
{\sl deformed} generators (in each context) in terms of the undeformed ones.

\newpage

{\bf Acknowledgments}  

We want to thank J.F. Cornwell for drawing our attention on 
ref. \cite{Zachos}. This 
work has been elaborated under a TOURNESOL financial support 
(which is also acknowledged by all of us).

\end{document}